%% file: IntroductoryLecture-final.tex
\documentclass[twoside,12pt, leqno]{article}
\usepackage{amsmath,amscd,amsthm,amssymb,amsxtra,latexsym,epsfig,epic,graphics}
\usepackage[matrix,arrow,curve]{xy}
\usepackage{graphicx}
\usepackage[breaklinks,bookmarksopen,bookmarksnumbered]{hyperref}
\hypersetup{colorlinks=true,backref=true,citecolor=blue,}
\input preamble.tex

\newcommand\<{{\langle}}
\renewcommand\>{{\rangle}}

\def\gM{{\mathfrak M}}
\def\gP{{\mathfrak P}}
\def\gW{{\mathfrak W}}
\def\AA{{\mathbb A}}

\def\PP{{\mathbb P}}
\def\GG{{\mathbb G}}

\def\QQ{{\mathbb Q}}
\def\FF{{\mathbb F}}

\def\Hilb{{\rm Hilb}}

\def\Spec{{{\rm Spec}\,}}

\numberwithin{equation}{section}

\usepackage{times}
\newdimen\x \x=12pt

\usepackage{color}

\usepackage[hang,flushmargin]{footmisc} 

\title{Syzygies, finite length modules, 
\\ and random curves}
\author{Christine Berkesch and 
Frank-Olaf Schreyer\footnote{
This article 
consists of extended notes from lectures by 
the second author at the Joint Introductory 
Workshop: Cluster Algebras and Commutative 
Algebra at MSRI during the Fall of 2012.  
This was supported by the National Science Foundation under Grant No. 0932078 000. 
}}


\begin{document}
\date{}
\maketitle
\begin{abstract}
We apply the theory of Gr\"obner bases to the 
computation of free resolutions over a 
polynomial ring, the defining equations of a 
canonically embedded curve, and the 
unirationality of the moduli space of curves of 
a fixed small genus.
\end{abstract}

\section*{Introduction}
\label{sec:intro}

While a great deal of modern commutative algebra 
and algebraic geometry has taken a 
nonconstructive form, the theory of Gr\"obner 
bases provides an algorithmic approach. 
Algorithms currently implemented in computer 
algebra systems, such as {\it Macaulay2} and 
{\sc Singular}, already exhibit the wide range 
of computational possibilities that arise from 
Gr\"obner bases~\cite{M2,Singular}. 

In these lectures, we focus on certain 
applications of Gr\"obner bases to syzygies and 
curves. 
In Section~\ref{sec:syzygy thm}, we use 
Gr\"obner bases to give an algorithmic proof of 
Hilbert's Syzygy Theorem, which bounds the length of a 
free resolution over a polynomial ring. 
In Section~\ref{sec:petri}, we prove Petri's 
theorem about the defining equations for 
canonical embeddings of curves. 
We turn in Section~\ref{sec:HRmods} to the 
Hartshorne--Rao module of a curve, showing by 
example how a module $M$ of finite length can be 
used to explicitly construct a curve whose 
Hartshorne--Rao module is $M$. 
Section~\ref{sec:unirational} then applies this 
construction to the study of the unirationality 
of the moduli space $\mathfrak M_g$ of curves of 
genus $g$. 

\newpage
\section{Hilbert's Syzygy Theorem}
\label{sec:syzygy thm}

Let $R:=\bK[x_1,\ldots,x_n]$ be a polynomial 
ring in $n$ variables over a field $\bK$. A 
\textbf{free resolution} of a finitely generated 
$R$-module $M$ is a complex of free modules 
$\cdots \to R^{\beta_2}\to R^{\beta_1}\to 
R^{\beta_0}$ such that the following is exact: 
\[
\cdots\to R^{\beta_2}\to R^{\beta_1}
\to R^{\beta_0}\to M\to 0. 
\]

\begin{theorem}[Hilbert's Syzygy Theorem]
\label{thm:syzygy-theorem}
\index{Hilbert!syzygy theorem}\index{theorem!syzygy}
Let $R=\PR$ be a polynomial ring in $n$ 
variables over a field $\bK$. Every finitely 
generated $R$-module $M$ has a finite free 
resolution of length at most $n$.
\end{theorem}

In this section, we give an algorithmic 
Gr\"obner basis proof of 
Theorem~\ref{thm:syzygy-theorem}, whose strategy 
is used in modern computer algebra systems like  
\textit{Macaulay2} and {\sc Singular} for syzygy 
computations. Gr\"obner bases were introduced by 
Gordan to provide a new proof of Hilbert's Basis 
Theorem~\cite{Gord99}. We believe that Gordan 
could have given the proof of 
Theorem~\ref{thm:syzygy-theorem} presented here.

\begin{definition}
\label{def:monomial order}
A \textbf{(global) monomial order} on $R$ is a 
total order $>$ on the set of monomials in $R$ 
such that the following two statements hold:
\begin{enumerate}
\vspace{-2mm}
\item
if $x^\alpha > x^\beta$, then 
$x^\gamma x^\alpha > x^\gamma x^\beta$ for all 
$\gamma\in\NN^n$, and
\vspace{-3mm}
\item
$x_i > 1$ for all $i$.
\end{enumerate}
\vspace{-2.5mm}
Given a global monomial order, the 
\textbf{leading term} of a nonzero polynomial 
$f = \sum_{\alpha} f_\alpha x^\alpha \in R$ is 
defined to be
\[
\ini(f) := f_\beta x^\beta,
\quad \text{where }
x^\beta := \max_\alpha \{ x^\alpha \mid f_\alpha\not=0\}.
\] 
For convenience, set $\ini(0) := 0$.
\end{definition}

\begin{theorem}[Division with Remainder]
\label{thm:det-div}
Let $>$ be a global monomial order on $R$, and 
let $\fr \in R$ be nonzero polynomials. For 
every $g\in R$, there exist uniquely determined 
$g_1,\ldots, g_r \in R$ and a \textbf{remainder} 
$h \in R$ such that the following hold.
\begin{itemize}
\vspace{-2.5mm}
\item[(1)]\quad
We have 
$g=g_1f_1+\cdots+g_rf_r + h$.
\vspace{-2.5mm}
\item[(2a)]\quad 
No term of $g_i\,\ini(f_i)$ is divisible by any 
$\ini(f_j)$ with $j<i$.
\vspace{-2.5mm}
\item[(2b)]\quad 
No term of $h$ is a multiple of $\ini(f_i)$ for 
any $i$.
\end{itemize}
\end{theorem}
\begin{proof} 
The result is obvious if $\fr$ are monomials, or 
more generally, if each $f_i$ has only a single 
nonzero term. Thus there is always a unique 
expression
\[
g = \sum_{i=1}^r g^{(1)}_i\,\ini(f_i) + h^{(1)}, 
\]
if we require that $g^{(1)}_1,\ldots, g^{(1)}_r$ 
and $h^{(1)}$ satisfy (2a) and (2b). By 
construction, the leading terms of the summands 
of the expression
\[
\sum_{i=1}^r g^{(1)}_i f_i+ h^{(1)}
\] 
are distinct, and the leading term in the 
difference 
$g^{(1)}
=g-(\sum_{i=1}^r g^{(1)}_i\,f_i + h^{(1)})$ 
cancels. Thus $\ini(g^{(1)}) < \ini(g)$, and we 
may apply induction on the number of summands in 
$g$.
\end{proof}

The remainder $h$ of the division of $g$ by 
$\fr$ depends on the order of $\fr$, since the 
partition of the monomials in $R$ given by (2a) 
and (2b) depends on this order. 
Even worse, it might not be the case that if 
$g \in \< \fr \>$, then $h=0$. 
A Gr\"obner basis is a system of generators for 
which this desirable property holds.

\begin{definition}
Let $I \subset R$ be an ideal. The 
\textbf{leading ideal} of $I$ (with respect to a 
given global monomial order) is 
\[
\ini(I) := \< \,\ini(f) \mid f \in I\, \>.
\]
A finite set $\fr$ of polynomials is a 
\textbf{Gr\"obner basis} when
\[
\< \ini (\< \fr \>) \>
=\< \, \ini (f_1),\ldots,\ini (f_r)\,\>.
\]
\end{definition}

Gordan's proof of Hilbert's basis theorem now 
follows from the easier statement that monomial 
ideals are finitely generated. In combinatorics, 
this result is called 
Dickson's Lemma~\cite{Dixon}.

If $\fr$ is a Gr\"obner basis, then by 
definition, a polynomial $g$ lies in $\< \fr \>$ 
if and only if the remainder $h$ under division 
of $g$ by $\fr$ is zero. In particular, in this 
case, the remainder does not depend on the order 
of $\fr$, and the monomials 
$x^\alpha \notin 
\< \ini(f_1),\ldots,\ini(f_r) \>$ represent a 
$\bK$-vector space basis of the quotient ring 
$R/\< \fr \>$, a fact known as 
Macaulay's theorem~\cite{Mac16}. For these 
reasons, it is desirable to have a Gr\"obner 
basis on hand. 

The algorithm that computes a Gr\"obner basis 
for an ideal is due to 
Buchberger~\cite{Buch,Buch2}. Usually, 
Buchberger's Criterion is formulated in terms of 
so-called S-pairs. In the treatment below, we do 
not use S-pairs; instead, we focus on the 
partition of the monomials of $R$ induced by 
$\ini(f_1), \ldots, \ini(f_r)$ via (2a) and (2b) 
of Theorem~\ref{thm:det-div1}.

Given polynomials $\fr$, consider the monomial 
ideals
\begin{align}
\label{eq:Mi}
M_i := \< \ini(f_1),\ldots,\ini(f_{i-1})\> : \ini(f_i) \quad \text{for \,$i=1,\ldots,r$}.
\end{align}
For each minimal generator $x^\alpha$ of an 
$M_i$, let $h^{(i,\alpha)}$ denote the remainder 
of $x^\alpha f_i$ divided by $\fr$ (in this 
order!).

\begin{theorem}[Buchberger's Criterion~\cite{Buch2}]
\label{thm:buchberger}
Let $\fr\in R$ be a collection of nonzero 
polynomials. Then $\fr$ form a Gr\"obner basis 
if and only if all of the remainders 
$h^{(i,\alpha)}$ are zero.
\end{theorem}

We will prove this result after a few more 
preliminaries. 

\begin{example}
\label{ex:generic3x5matrix} 
Consider the ideal generated by the $3\times 3$ 
minors of the matrix
\[
\begin{pmatrix}
x_1 & x_2 & x_3 & x_4 & x_5 \cr
y_1 & y_2 & y_3 & y_4 & y_5 \cr
z_1 & z_2 & z_3 & z_4 & z_5 \cr
\end{pmatrix}. 
\]
Using the lexicographic order on 
$\bK[x_1,\ldots,z_5]$, the leading terms of the 
maximal minors of this matrix and the minimal 
generators of the corresponding monomial ideals 
$M_i$ are listed in the following table. 
\[
\begin{array}{|c|l|}
\hline
\quad x_1y_2z_3\quad &\quad M_1=0 \quad  \cr
\quad x_1y_2z_4\quad  & \quad M_2=\< z_3\>\quad\cr 
\quad x_1y_3z_4\quad  & \quad M_3=\< y_2\>\quad\cr
\quad x_2y_3z_4\quad  & \quad M_4=\< x_1\>\quad\cr
\quad x_1y_2z_5\quad  & \quad M_5=\< z_3,z_4\>\quad\cr 
\quad x_1y_3z_5\quad  & \quad M_6=\< y_2,z_4\>\quad\cr
\quad x_2y_3z_5\quad  & \quad M_7=\< x_1,z_4\>\quad\cr
\quad x_1y_4z_5\quad  & \quad M_8=\< y_2,y_3\>\quad\cr
\quad x_2y_4z_5\quad  & \quad M_9=\< x_1,y_3\>\quad\cr
\quad x_3y_4z_5\quad  & \quad M_{10}=\< x_1,x_2\>\quad\cr
\hline
\end{array}
\]
Note that only 15 of the possible 
${10 \choose 2}=45$ S-pairs are needed to test 
Buchberger's Criterion, as stated in 
Theorem~\ref{thm:buchberger}. 
\end{example}

\begin{exercise}
Show that the maximal minors of the matrix in 
Example~\ref{ex:generic3x5matrix} form a 
Gr\"obner basis by using the Laplace expansions 
of suitable $4\times 4$ matrices.
\end{exercise}

In order to prove 
Theorems~\ref{thm:syzygy-theorem} 
and~\ref{thm:buchberger}, we now extend the 
notion of a monomial order to vectors of 
polynomials. 

\begin{definition}
\label{def:order free}
A monomial of a free module $R^r$ with basis 
$e_1,\ldots,e_r$ is an expression 
$x^\alpha e_i$. A 
\textbf{(global) monomial order} on $R^r$ is a 
total order of the monomials of $R^r$ such that 
the following two statements hold:
\begin{enumerate}
\vspace{-2mm}
\item 
if $x^\alpha e_i > x^\beta e_j$, then 
$x^\gamma x^\alpha e_i > x^\gamma x^\beta e_j$  
for all $i,j$ and $\gamma\in\NN^n$, and
\vspace{-3mm}
\item
$x^\alpha e_i >e_i$ for all $i$ and 
$\alpha\not=0$.
\end{enumerate}
\end{definition}

Usually, it is also the case that 
$x^\alpha e_i > x^\beta e_i$ if and only if 
$x^\alpha e_j> x^\beta e_j$,
i.e., the order on the monomials in the 
components induce a single monomial order on 
$R$. 

Thanks to Definition~\ref{def:order free}, we 
may now speak of the leading term of a vector of 
polynomials. In this situation, the division 
theorem still holds.

\begin{theorem}[Division with Remainder for Vectors of Polynomials]
\label{thm:det-div1}
Let $>$ be a global monomial order on $R^{r_0}$, 
and let $F_1,\dots,F_r \in R^{r_0} $ be nonzero 
polynomial vectors. For every $G\in R^{r_0}$, 
there exist uniquely determined 
$g_1,\ldots,g_r \in R$ and a \textbf{remainder} 
$H \in R^{r_0}$ such that the following hold.
\begin{itemize}
\vspace{-3mm}
\item[(1)]\quad
We have $G=g_1F_1+\cdots+g_rF_r + H$.
\vspace{-3mm}
\item[(2a)]\quad 
No term of $g_i\, \ini(F_i)$ is a multiple of an 
$\ini(F_j)$ with $j< i$.
\vspace{-3mm}
\item[(2b)]\quad 
No term of $H$ is a multiple of $\ini(F_i)$ for 
any $i$.
\qed
\end{itemize}
\end{theorem}

\begin{definition}
Generalizing the earlier definition, given a 
global monomial order on $R^r$, the 
\textbf{leading term} of a nonzero vector of 
polynomials $F = (\fr)$ is defined to be the 
monomial 
\[
\ini(F) := f_{\beta_i} x^{\beta_i} e_i, 
\quad \text{where $x^{\beta_i} = \max_{\alpha_i}
\{x^{\alpha_i}\mid f_{\alpha_i}x^{\alpha_i} 
\text{ is a nonzero term of } f_i \}$.}
\] 
A finite set $F_1,\ldots,F_s$ of vectors of 
polynomials in $R^r$ is a 
\textbf{Gr\"obner basis} when 
\[
\< \ini(\<F_1,\ldots,F_s\>) \> 
= \<\,\ini(F_1),\ldots,\ini(F_s)\,\>. 
\]
\end{definition}

We are now prepared to prove 
Theorem~\ref{thm:buchberger}, followed by a 
corollary. 

\begin{proof}[Proof of Theorem~\ref{thm:buchberger} (Buchberger's Criterion)]
The forward direction follows by definition. For 
the converse, assume that all remainders 
$h^{(i,\alpha)}$ vanish. Then for each minimal 
generator $x^\alpha$ in an $M_i$, there is an 
expression 
\begin{equation}
\label{eq:rel-1}
x^\alpha f_i = 
g_1^{(i,\alpha)}f_1 +\cdots+ g_r^{(i,\alpha)}f_r
\end{equation}
such that no term of $g^{(i,\alpha)}_j\ini(f_j)$ 
is divisible by an $\ini(f_k)$ for every $k<j$, 
by condition (2a) of Theorem~\ref{thm:det-div}.
(Of course, for a suitable $j<i$, one of the 
terms of $g^{(i,\alpha)}_j\ini(f_j)$ coincides 
with $x^\alpha \ini(f_i)$. This is the second 
term in the usual S-pair description of 
Buchberger's Criterion.) 
Now let $e_1,\ldots, e_r \in R^r$ denote the 
basis of the free module, and let 
$\varphi:R^r\to R$ be defined by 
$e_i \mapsto f_i$. 
Then by~\eqref{eq:rel-1}, elements of the form 
\begin{align}
\label{eq:Gialpha}
G^{(i,\alpha)} := 
-g^{(i,\alpha)}_1e_1+\cdots+(x^\alpha-g^{(i,\alpha)}_i)e_i+\cdots + &(-g_r^{(i,\alpha)})e_r 
\end{align}
are in the kernel of $\phi$. In other words, the 
$G^{(i,\alpha)}$'s are syzygies between $\fr$. 

We now proceed with a division with remainder in 
the free module $R^r$, using the induced 
monomial order $>_1$ on $R^r$ defined by
\begin{align}
\label{eq:order free}
\begin{array}{ll} x^\alpha e_i >_1 x^\beta e_j 
\iff
& x^\alpha \ini(f_i) 
> x^\beta \ini
(f_j) \;\hbox{ or }\;\cr & x^\alpha \ini(f_i) 
= x^\beta \ini(f_j) 
\;\hbox{ (up to a scalar)}\;\hbox{ with }\; 
i > j. \cr \end{array}
\end{align}
With respect to this order, 
\[
\ini(G^{(i,\alpha)})= x^\alpha e_i
\]
because the term $cx^\beta \ini(f_j)$ that 
cancels against $x^\alpha \ini(f_i)$ 
in~\eqref{eq:rel-1} satisfies $j<i$, and all 
other terms of any $g^{(i,\alpha)}_k \ini(f_k)$ 
are smaller.

Now consider an arbitrary element
\[
g = a_1f_1+\cdots+a_rf_r \in \< \fr \>.
\]
We must show that 
$\ini(g) \in \< \ini(f_1),\ldots,\ini(f_r) \>$.
Let $g_1e_1+\cdots+g_re_r$ be the remainder of 
$a_1e_1+\cdots +a_re_r$ divided by the 
collection of $G^{(i,\alpha)}$ vectors. Then 
\[
g=a_1f_1+\cdots+a_rf_r =g_1f_1+\cdots+g_rf_r 
\]
because the $G^{(i,\alpha)}$ are syzygies,
and $g_1,\ldots, g_r$ satisfy (2a) of 
Theorem~\ref{thm:det-div} when $a_1,\ldots,a_n$ 
are divided by $\fr$, by the definition of the 
$M_i$ in~\eqref{eq:Mi}.
Therefore, the nonzero initial terms
\[
\ini(g_jf_j)=\ini(g_j)\ini(f_j)
\]
are distinct and no cancellation can occur among 
them. The proof is now complete because 
\[
\ini(g) := \max_{j} \{ \ini(g_j)\ini(f_j) \} 
\in \< \ini(f_1),\ldots,\ini(f_r) \>.
\hfill\qedhere
\]
\end{proof}

\begin{corollary}[Schreyer~\cite{Sch0}]
\label{cor:Sch80}
If $F_1,\ldots,F_{r_1} \in R^{r_0}$ are a 
Gr\"obner basis, then the $G^{(i,\alpha)}$ 
of~\eqref{eq:Gialpha} form a Gr\"obner basis of
$\ker(\varphi_1\colon R^{r_1} \to R^{r_0})$ with 
respect to the induced monomial order $>_1$ 
defined in~\eqref{eq:order free}. In particular, 
$F_1,\ldots,F_{r_1}$ generate the kernel of 
$\varphi_1$.
\end{corollary}
\begin{proof} As mentioned in the proof of 
Theorem~\ref{thm:buchberger}, the coefficients 
$g_1, \ldots ,g_r$ of a remainder 
$g_1e_1+\cdots+g_re_r$ resulting from division 
by the $G^{(i,\alpha)}$ satisfy condition (2a) 
of Theorem~\ref{thm:det-div} when divided by 
$\fr$. Hence, no cancellation can occur in the 
sum $g_1 \ini(f_1)+\cdots+g_r \ini(f_r)$, and 
$g_1f_1+\cdots g_rf_r=0$ only if 
$g_1=\ldots=g_r=0$. Therefore, the collection of 
$\ini(G^{(i,\alpha)})$ generate 
the leading term ideal $\ini(\ker \varphi_1)$.
\end{proof}

We have reached the goal of this section, an 
algorithmic proof of 
Theorem~\ref{thm:syzygy-theorem}. 

\begin{proof}[Proof of Theorem~\ref{thm:syzygy-theorem} (Hilbert's Syzygy Theorem)]
Let $M$ be a finitely generated $R$-module with 
presentation 
\[
R^{r}\overset{\varphi}\lra R^{r_0}\lra M\lra 0. 
\]
Regard $\varphi$ as a matrix and, thus, its 
columns as a set of generators for 
$\im \varphi$. Starting from these generators, 
compute a minimal Gr\"obner basis 
$F_1, \ldots, F_{r_1}$ for $\im\varphi$ with 
respect to some global monomial $>_0$ order on 
$R^{r_0}$. Now consider the induced monomial 
order $>_1$ on $R^{r_1}$, and let 
$G^{(i,\alpha)}\in R^{r_1}$ denote the syzygies 
obtained by applying Buchberger's Criterion to 
$F_1,\ldots,F_{r_1}$. By 
Corollary~\ref{cor:Sch80}, the $G^{(i,\alpha)}$ 
form a Gr\"obner basis for the kernel of the map 
${\varphi_1}\colon R^{r_1} \rightarrow R^{r_0}$, 
so we may now repeat this process. 

Let $\ell$ be the maximal $k$ such that the 
variable $x_k$ occurs in some leading term 
$\ini(F_j)$. Sort $F_1,\ldots,F_{r_1}$ so that 
whenever $j < i$, the exponent of $x_\ell$ in 
$\ini(F_j)$ is less than or equal to the 
exponent of $x_\ell$ in $\ini(F_i)$. In this 
way, none of the variables $x_\ell,\ldots, x_n$ 
will occur in a leading term 
$\ini(G^{(i,\alpha)})$. Thus the process will 
terminate after at most $n$ steps.
\end{proof}

Note that there are a number of choices allowed 
in the algorithm in the proof of 
Theorem~\ref{thm:syzygy-theorem}. In particular, we may 
order each set of Gr\"obner basis elements as we 
see fit. 

\begin{example}
\label{ex:5gon}  
Consider the ideal
$I=\< f_1 , \dots ,f_5\>\subset R=\bK[w,x,y,z]$ 
generated by the polynomials
\[
f_1=w^2-xz,\,\,f_2=wx-yz,\,\,f_3=x^2-wy,\,\,f_4=xy-z^2,\,\,f_5=y^2-wz.
\]
To compute a finite free resolution of $M=R/I$ 
using the method of the proof of 
Theorem~\ref{thm:syzygy-theorem}, we use the 
degree reverse lexicographic order on $R$. The 
algorithm successively produces three syzygy 
matrices $\varphi_1$, $\varphi_2$, and 
$\varphi_3$, which we present in a compact way 
as follows. 
\[
\begin{matrix}
{\bf w^2}-xz &\vline& -x & y & 0  &-z & 0 & -y^2+wz \cr
{\bf wx}-yz&\vline&{\bf w} &-x &-y & 0 & z & z^2     \cr
{\bf x^2}-wy&\vline&-z &{\bf w} & 0 &-y & 0 &  0      \cr
{\bf xy}-z^2&\vline& 0  & 0 &{\bf w} &{\bf x }&-y & yz      \cr
 {\bf y^2}-wz&\vline& 0 & 0 &-z &-w &{\bf x} &{\bf w^2} \cr
\hline
&\vline& 0 & y&-x&{\bf w}&-z& 1 \cr
&\vline&-y^2+wz & z^2 & -wy & yz & -w^2 &{\bf x}
\end{matrix}
\]
All initial terms are printed in bold. The first 
column of this table is the transpose of the 
matrix $\varphi_1$. It contains the original 
generators for $I$ which, as Buchberger's 
Criterion shows, already form a Gr\"obner basis 
for $I$. The syzygy matrix $\varphi_2$  
resulting from the algorithm is the 
$5\times 6$ matrix in the middle of our table. 
Note that, for instance, $M_4 = \< w,x\>$ can be 
read from the $4$th row of $\varphi_2$. 

By Corollary~\ref{cor:Sch80}, we know that the 
columns of $\varphi_2$ form a Gr\"obner basis 
for $\ker(\varphi_1)$ with respect to the 
induced monomial order on $R^5$. Buchberger's 
Criterion applied to these Gr\"obner basis 
elements yields a $6 \times 2$ syzygy matrix 
$\varphi_3$, whose transpose is printed in the 
two bottom rows of the table above. Note that 
there are no syzygies on the two columns of 
$\varphi_3$ because the initial terms of these 
vectors lie with different basis vectors. 

To summarize, we obtain a free resolution of the 
form
\[ 
0 \longrightarrow R^2 \overset{\varphi_3}\longrightarrow R^6 \overset{\varphi_2}
\longrightarrow R^5 \overset{\varphi_1}\longrightarrow R \longrightarrow
R/I \longrightarrow 0.
\]
Observe that, in general, once we have the 
initial terms of a Gr\"obner basis for $I$, 
we can easily compute the initial terms of the 
Gr\"obner bases for all syzygy modules, that is, 
all bold face entries of our table. This gives 
us an idea on the amount of computation that 
will be needed to obtain the full free 
resolution. 
\end{example}

If the polynomial ring is graded, say 
$R=S=\bK[x_0,\ldots,x_n]$ is the homogeneous 
coordinate ring of $\PP^n$, and $M$ is a 
finitely generated graded $S$-module, then the 
resolution computed through the proof of 
Theorem~\ref{thm:syzygy-theorem} is homogeneous 
as well. However, this resolution is typically 
not minimal. In Example~\ref{ex:5gon}, the last 
column of $\varphi_2$ is in the span of the 
previous columns, as can be seen from the first 
row of $\varphi_3^t$.

\begin{example} 
Recall that in 
Example~\ref{ex:generic3x5matrix}, we considered 
the ideal $I$ of $3\times 3$ minors of a generic 
$3\times 5$ matrix over $S=\bK[x_1,\ldots,z_5]$ 
with the standard grading. 
The algorithm in the proof of 
Theorem~\ref{def:monomial order} produces a 
resolution of $S/I$ of the form 
\[
0 \longrightarrow S(-5)^6 \longrightarrow S(-4)^{15} \longrightarrow S(-3)^{10} 
\longrightarrow S \longrightarrow S/I \longrightarrow 0
\]
because $I$ is generated by $10$ Gr\"obner basis 
elements, there are altogether $15$ minimal 
generators of the $M_i$ ideals, and $6$ of the 
monomial ideal $M_i$ have $2$ generators. In 
this case, the resolution is minimal for degree 
reasons.
\end{example}

\begin{exercise} 
Let $I$ be a Borel-fixed monomial ideal. Prove 
that in this case, the algorithm in the proof of 
Theorem~\ref{thm:syzygy-theorem} produces a 
minimal free resolution of $I$. Compute the 
differentials explicitly and compare your result 
with the complex of S. Eliahou and Kervaire 
in~\cite{EK} (see also~\cite{PS08}).
\end{exercise}

\section{Petri's Theorem}
\label{sec:petri}

One of the first theoretical applications of 
Gr\"obner bases is Petri's analysis of the 
generators of the homogeneous ideal of a 
canonically embedded curve. Petri was the last 
student of Max Noether, and he acknowledges  
help from Emmy Noether in his thesis. As Emmy 
Noether was a student of Gordan, it is quite 
possible that Petri became aware of the concept 
of Gr\"obner bases through his communication 
with her, but we do not know if this was the case.

Let $C$ be a smooth projective curve of genus 
$g$ over $\CC$. Let 
$\omega_1, \ldots,\omega_g \in H^0(C, \omega_C)$ 
be a basis of the space of holomorphic 
differential forms on $C$ and consider the 
canonical map 
\[
\iota: C \to \PP^{g-1}
\quad \text{given by}\quad
p\mapsto [\omega_1( p): \cdots: \omega_g( p)].
\]
The map $\iota$ is an embedding unless $C$ is 
hyperelliptic. We will assume that $C$ is not 
hyperelliptic. 
Let $S:=\CC[x_1,\ldots,x_g]$ be the homogeneous 
coordinate ring of $\PP^{g-1}$, and let $I_C \subset S$ be the homogeneous ideal of $C$.

\begin{theorem}[Petri's Theorem~\cite{Pet}]
\label{thm:Petri} 
The homogeneous ideal of a canonically embedded 
curve is generated by quadrics unless
\begin{itemize} 
\vspace{-2.5mm}
\item $C$ is trigonal (i.e., there is a 3:1 
holomorphic map $C \to \PP^1$) or
\vspace{-2.5mm}
\item $C$ is isomorphic to a smooth plane 
quintic. In this case, $g=6$. 
\end{itemize}
\end{theorem}

Petri's Theorem received much attention through the 
work of Mark Green \cite{Gre}, who formulated a conjectural 
generalization to higher syzygies of canonical 
curves in terms of the Clifford index. We will 
not report here on the impressive progress made 
on this conjecture in the last two decades, but 
refer instead 
to~\cite{AF,AN,AV,GL,HR,Muk,Sch1,Sch2,Sch3,Vo1,Vo2, Vo3} for 
further reading. 

In the cases of the exceptions in 
Theorem~\ref{thm:Petri}, also Babbage~\cite{Ba}  observed 
 that the ideal cannot be generated 
by quadrics alone. If 
$D := p_1+\cdots+p_d$ is a divisor of degree $d$ 
on $C$, then the linear system $|\omega_C(-D)|$ 
is cut out by hyperplanes through the span 
$\overline D$ of the points 
$p_i \in C \subset \PP^{g-1}$. Thus 
Riemann--Roch implies that 
\[
h^0(C,\sO_C(D))= d+1- g +  \codim \overline D = d- \dim \overline D.
\]
Hence the three points of a trigonal divisor 
span only a line, and by B\'ezout's theorem, we 
need cubic generators in the generating set of 
its vanishing ideal. 

Similarly, in the second exceptional case, the 5 
points of a $g^2_5$ are contained in a unique 
conic in the plane they span, and quadrics 
alone do not cut out the curve. 

The first step of Petri's analysis builds upon a 
proof by Max Noether.

\begin{theorem}[M.~Noether~\cite{Noe}]
\label{thm:Noether}
A non-hyperelliptic canonical curve 
$C \subset \PP^{g-1}$ is projectively normal, 
i.e., the maps 
$H^0(\PP^{g-1},\sO(n))\to H^0(C,\omega_C^{\otimes n})$ 
are surjective for every $n$.
\end{theorem}

\begin{proof} Max Noether's proof is a clever 
application of the base point free pencil trick. 
This is a method which, according to Mumford, 
Zariski taught to all of his students. Let $|D|$ 
be a base point free pencil on a curve, and let 
$\sL$ be a further line bundle on $C$. Then the 
Koszul complex
\[
0 \to \Lambda^2 H^0(\sO_C(D)) \tensor \sL(-D) 
\to H^0(\sO_C(D)) \tensor \sL \to \sL(D) \to 0
\] 
is an exact sequence. To see this, note that 
locally, at least one section of the line bundle 
$\sO_C(D)$ does not vanish. Thus the kernel of 
the multiplication map
\[
H^0(\sO_C(D)) \tensor H^0(\sL)  \to H^0(\sL(D))
\] 
is isomorphic with $H^0(\sL(-D))$. (Note 
$\Lambda^2 H^0(\sO_C(D)) \cong \CC$, as 
$h^0(\sO_C(D))=2$.)
 
Consider $p_1, \ldots, p_g$ general points on 
$C$ and the divisor $D=p_1+\cdots+p_{g-2}$ built 
from the first $g-2$ points. Then the images of 
these points span $\PP^{g-1}$ and the span of 
any subset of less than $g-1$ points intersects 
the curve in no further points. Choose a basis 
$\omega_1, \ldots, \omega_g \in H^0(\omega_C)$ 
that is, up to scalars, dual to these points, 
i.e., $\omega_i(p_j) = 0$ for $i\not= j$ and 
$\omega_i(p_i)\neq 0$. Then $|\omega_C(-D)|$ is 
a base point free pencil spanned by 
$\omega_{g-1}, \omega_g$. If we apply the base 
point free pencil trick to this pencil and 
$\sL=\omega_C$, then we obtain the sequence
\[ 
0 \to 
\Lambda^2 H^0\omega_C(-D) \tensor H^0\sO_C(D) 
\to H^0\omega_C(-D) \tensor H^0\omega_C 
\longrightarrow^{\hspace{-4.5mm}\mu
\hspace{2mm}} H^0\omega_C^{\tensor 2}(-D),
\]
and the image of 
\begin{equation}
\label{eq:mult}
\mu\colon 
H^0(\omega_C(-D)) \tensor H^0(\omega_C) \to 
H^0(\omega_C^{\tensor 2}(-D))
\end{equation} 
is $2g-1$ dimensional because 
$h^0(\omega_C(-D))=2$ and $h^0(\sO_C(D))=1$. 
Thus $\mu$ in~\eqref{eq:mult} is 
surjective, since 
$h^0(\omega_C^{\tensor 2}(-D))=2g-1$ holds by 
Riemann--Roch. On the other hand,
\[
\omega_1^{\tensor 2}, \ldots, 
\omega_{g-2}^{\tensor 2} \in 
H^0(\omega_C^{\tensor 2})
\] 
represent linearly independent elements of
$H^0(\omega_C^{\tensor 2})/H^0(\omega_C^{\tensor 2}(-D))$, 
hence represent a basis, and the map 
$H^0(\omega_C) \tensor H^0(\omega_C) 
\to H^0(\omega_C^{\tensor 2})$
is surjective as well. This proves quadratic 
normality. 
 
The surjectivity of the multiplication maps
\[
H^0(\omega_C^{\tensor n-1})\tensor H^0(\omega_C) 
\to H^0(\omega_C^{\tensor n})
\]
for $n\ge 3$ is similar, but easier: 
$\omega_1^{\tensor n},\ldots,
\omega_{g-2}^{\tensor n}\in 
H^0(\omega_C^{\tensor n})$ are linearly 
independent modulo the codimension $g-2$ 
subspace $H^0(\omega_C^{\tensor n}(-D))$, 
and the map
\[ 
H^0(\omega_C^{\tensor n-1}) \tensor 
H^0(\omega_C(-D)) 
\to H^0(\omega_C^{\tensor n}(-D))
\] 
is surjective simply because 
$H^1(\omega_C^{\tensor n-2}(D))=0$ for $n\ge 3$.
\end{proof}

\begin{corollary} 
The Hilbert function of the coordinate ring of a 
canonical curve takes values
\vspace{-3mm}
\[
\hspace{2.9cm}
\dim (S/I_C)_n = \begin{cases}
1 & \hbox{ if } n= 0\\
g & \hbox{ if } n= 1\\
(2n-1)(g-1) &\hbox{ if } n\ge 2.
\hspace{2.28cm}\qed
\end{cases} 
\] 
\end{corollary} 

\begin{proof}[Proof of Theorem~\ref{thm:Petri} (Petri's Theorem)]
Petri's analysis begins with the map $\mu$ 
in~\eqref{eq:mult} above. Choose homogeneous 
coordinates $x_1,\ldots, x_g$ such that 
$x_i \mapsto \omega_i$. Since 
$\omega_i\tensor \omega_j \in 
H^0(\omega_C^{\tensor 2}(-D))$ for 
$1 \le i <j \le g-2$, we find the polynomials 
\begin{equation}\label{eq:fij}
f_{ij} := x_ix_j -\sum_{r=1}^{g-2} 
a_{ij}^r x_r - b_{ij} \in I_C, 
\end{equation}
where the $a^r_{ij}$ and $b_{ij}$ are linear and  
quadratic, respectively, in $\CC[x_{g-1},x_g]$. 
We may choose a monomial order such that 
$\ini(f_{ij})=x_ix_j$. Since 
${g-2 \choose 2}={g+1 \choose 2} - (3g-3)$, 
these quadrics span $(I_C)_2$. On the other 
hand, they do not form a Gr\"obner basis for 
$I_C$ because the $(g-2){n+1 \choose 2}+(n+1)$ 
monomials $x_i^kx_{g-1}^\ell x_g^m$ with 
$i=1, \ldots, g-2$ and  $k+\ell+m=n$ represent a 
basis for 
$(S/\< x_ix_j\mid 1\leq i<j\leq g-2 \>)_n$, 
which is still larger. We therefore need $g-3$ 
further cubic Gr\"obner basis elements. To find 
these, Petri considers the base point free 
pencil trick applied to $|\omega_C(-D)|$ and 
$\sL=\omega_C^{\tensor 2}(-D)$. The cokernel of 
the map
\begin{equation}\label{eq:mult3}
H^0(\omega_C(-D)) \tensor 
H^0(\omega_C^{\tensor 2}(-D)) 
\to H^0(\omega_C^{\tensor 3}(-2D))
\end{equation}
has dimension $h^1(\omega_C)=1$. To find the 
missing element in 
$H^0(\omega_C^{\tensor 3}(-2D))$, Petri 
considers the linear form 
$\alpha_i=\alpha_i(x_{g-1},x_g)$ in the pencil 
spanned by $x_{g-1},x_g$ that defines a tangent 
hyperplane to $C$ at $p_i$. Then 
$\alpha_i\omega_i^{\tensor 2} \in 
H^0(\omega_C^{\tensor 3}(-2D))$ because 
$\omega_i^{\tensor 2}$ vanishes quadratically at 
all points $p_j \not=p_i$, while $\alpha_i$ 
vanishes doubly at $p_i$. Not all of these 
elements can be contained in the image 
of~\eqref{eq:mult3}, since otherwise we would 
find $g-2$ further cubic Gr\"obner basis 
elements of type
\[
\alpha_i x_i^2 + \hbox{ lower order terms,}
\]
where a lower order term is a term that is at 
most linear in $x_1,\ldots,x_{g-2}$. As  this is 
too many, at least one of the 
$\alpha_i\omega_i^{\tensor 2}$ spans the 
cokernel of the map~\eqref{eq:mult3}. 

We now argue by uniform position. Since $C$ is 
irreducible, the behavior of 
$\alpha_i\omega_i^{\tensor 2}$ with respect to 
spanning of the cokernel is the same for any 
general choice of points $p_1, \ldots,p_g$. So 
for general choices, each of these elements span 
the cokernel, and after adjusting scalars, we 
find that 
\begin{equation}
\label{eq:Gkl}
G_{k\ell} := \alpha_k x_k^2 -\alpha_\ell 
x_\ell^2 + \hbox{ lower order terms }
\end{equation}
are in $I_C$. Note that 
$G_{k\ell}= - G_{\ell k}$ and 
$G_{k\ell}+G_{\ell m} = G_{km}$. So this gives 
only $g-3$ further equations with leading terms
$x_k^2x_{g-1}$ for $k=1,\ldots,g-3$ up to a 
scalar. The last Gr\"obner basis element is a 
quartic $H$ with leading term 
$\ini(H)=x_{g-2}^3x_{g-1}$, which we can obtain 
as a  remainder of the Buchberger test 
applied to $x_{g-2}G_{k,g-2}$. There are no 
further Gr\"obner basis elements, because the 
quotient $S/J$ of $S$ by
\[
J:=\< x_ix_j, x_k^2x_{g-1}, x_{g-2}^3x_{g-1} 
\mid 1 \le i< j \le g-2, 1 \le k \le g-3 \>
\] 
has the same Hilbert function as $S/I_C$. Hence 
$\ini(I_C) = J$. 

We now apply Buchberger's test to $x_kf_{ij}$ 
for a triple of distinct indices 
$1\le i,j,k\le g-2$. Division with remainder 
yields a syzygy
\begin{equation}\label{eq:Petrisyz}
x_kf_{ij} - x_jf_{ik} + \sum_{r\not=k} 
a_{ij}^r f_{rk} - \sum_{r \not= j} 
a_{ik}^r f_{rj} + \rho_{ijk} G_{kj} =0
\end{equation}
for a suitable coefficient $\rho_{ijk} \in \CC$. 
(Moreover, comparing coefficients, we find that 
$a_{ij}^k = \rho_{ijk}\alpha_k$ holds. In 
particular, Petri's coefficients $\rho_{ijk}$ 
are symmetric in $i,j,k$, since $a_{ij}^k$ is 
symmetric in $i,j$.) Since $C$ is irreducible, 
we have that for a general choice of 
$p_1, \ldots, p_g$, either all coefficients 
$\rho_{ijk} \not=0$ or all $\rho_{ijk}=0$. In 
the first case, the cubics lie in the ideal 
generated by the quadrics. 

In the second case, the $f_{ij}$ are a Gr\"obner 
basis by themselves. Thus the zero locus 
$V(f_{ij}| 1 \le i<j \le g-2)$ of the quadrics 
$f_{ij}$ define an ideal of a scheme $X$ of 
dimension 2 and degree $g-2$. Since $C$ is 
irreducible and non-degenerate, the surface $X$ 
is irreducible and non-degenerate as well. Thus 
$X \subset \PP^{g-2}$ is a surface of minimal 
degree. These were classified by Bertini, see, for instance,~\cite{EH85}.
 Either $X$ is a rational normal surface 
scroll, or $X$ is isomorphic to the Veronese 
surface $\PP^2 \hookrightarrow \PP^5$. In the 
case of a scroll, the ruling on $X$ cuts out a 
$g^1_3$ on $C$ by Riemann--Roch. In the case of 
the Veronese surface, the preimage of $C$ in 
$\PP^2$ is a plane quintic.  
\end{proof}

Perhaps the most surprising part of Petri's 
theorem is this: either $I_C$ is generated by 
quadrics or there are precisely $g-3$ minimal 
cubic generators. It is a consequence of the 
irreducibility of $C$ that no value in between 0 
and $g-3$ is possible for the number of cubic 
generators. If we drop the assumption of 
irreducibility, then there are canonical curves 
with $1,\ldots,g-5$ or $g-3$ cubic generators. 
For example, if we take a stable curve 
$C=C_1\cup C_2$ with two smooth components of 
genus $g_i \ge 1$ intersecting in three points, 
so that $C$ has genus $g=g_1+g_2+2$, then the 
dualizing sheaf $\omega_C$ is very ample and the 
three intersection points lie on a line by the 
residue theorem. For general curves $C_1$ and 
$C_2$ of genus $g_i \ge 3$ for $i\in\{1,2\}$, 
the ideal $I_C$ has precisely one cubic 
generator, see~\cite{Sch2}. However, we could 
not find such an example with precisely $g-4$ 
generators. For genus $g=5$, one cubic generator 
is excluded by the structure theorem of 
Buchsbaum--Eisenbud, and obstructions for larger 
$g$ are unclear to us. 

\begin{conjecture} 
\label{conj:Gor mingens}
Let $A=S/I$ be a graded artinian Gorenstein 
algebra with Hilbert function $\{1,g-2,g-2,1\}$. 
Then $I$ has $0,1,\ldots,g-5$ or $g-3$ cubic 
minimal generators.
\end{conjecture} 

The veracity of this conjecture would imply the 
corresponding statement for reducible canonical 
curves because the artinian reduction 
$A:=S/(I_C+\< \ell_1,\ell_2 \>)$ of $S/I_C$, 
for general linear forms $\ell_1,\ell_2$, 
has Hilbert function $\{1,g-2,g-2,1\}$.

Petri's analysis has been treated by 
Mumford~\cite{Mum}, as well 
as~\cite{ACGH,S-D,Sho}. From our point of view, 
Gr\"obner bases and the use of uniform position 
simplify and clarify the treatment quite a bit. 
Mumford remarks in~\cite{Mum} that we now have 
seen all curves at least once, following a claim 
of Petri~\cite{Pet}. We disagree with him in 
this point. If we introduce indeterminates for 
all of the coefficients in Petri's equations, 
then the scheme defined by the condition on the 
coefficients that $f_{ij}, G_{kl}$, and $H$ form 
a Gr\"obner basis can have many 
components~\cite{Sch2,Lit}. It is not clear to 
us how to find the component corresponding to 
smooth curves, much less how to find 
closed points on this component. 

\section{Finite length modules and space curves}
\label{sec:HRmods}

In the remaining part of these lectures, we  
report on how to find all curves in a Zariski 
open subset of the moduli space $\gM_g$ of 
curves of genus $g$ for small $g$. In 
Section~\ref{sec:unirational}, we report on the 
known \emph{unirationality} results for these 
moduli spaces. But first, we must discuss a 
method to explicitly construct space curves. 

In this section, a \emph{space curve} 
$C \subset \PP^3$ will be a Cohen--Macaulay 
subscheme of pure dimension $1$; in particular, 
$C$ has no embedded points. We denote by $\sI_C$ 
the ideal sheaf of $C$ and by 
$I_C = \sum_{n\in \ZZ} H^0(\PP^3, \sI_C(n))$ the 
homogeneous ideal of $C$. 
The goal of this section is to construct a curve 
$C$ of genus $g$ and degree $d$. To do so, we 
will use work of Rao, who showed that the 
construction of $C$ is equivalent to the 
creation of its Hartshorne--Rao module (see 
Theorem~\ref{thm:rao}). 

\begin{definition} 
\label{def:HRmod}
The \textbf{Hartshorne--Rao module} of $C$ is 
the finite length module
\[ 
M = M_C := \sum_{n \in \ZZ} H^1(\PP^3, \sI_C(n))
\subset \sum_{n\in \ZZ} H^0(\PP^3, \sO(n)) 
\cong S:=\bK[x_0,..,x_3].
\]
\end{definition}
The Hartshorne--Rao module measures the 
deviation of $C$ from being projectively normal. 
Furthermore, $M_C$ plays an important role in 
liaison theory of curves in $\PP^3$, which we 
briefly recall now. 

Let $S:=\bK[x_0,\ldots,x_3]$ and $S_C:=S/I_C$ 
denote the homogeneous coordinate ring of 
$\PP^3$ and $C \subset \PP^3$, respectively.
By the Auslander--Buchsbaum--Serre 
formula~\cite[Theorem 19.9]{Eis}, $S_C$ has 
projective dimension $\pd_S S_C \le 3$.
Thus its minimal free resolution has the form
\[ 
0 \leftarrow S_C \leftarrow S \leftarrow F_1  
\leftarrow F_2  \leftarrow F_3  \leftarrow 0, 
\]
with free graded modules 
$F_i=\oplus S(-j)^{\beta_{ij}}$.
By the same formula in the local case, we see 
that the sheafified 
$\sG:=\ker(\widetilde F_1 \to \sO_{\PP^3})$ 
is always a vector bundle, and
\begin{equation}
\label{OCres}
0\leftarrow\sO_C \leftarrow \sO_{\PP^3} 
\leftarrow 
\bigoplus_j \sO_{\PP^3}(-j)^{\beta_{1j}} 
\leftarrow \sG \leftarrow 0
\end{equation}
is a resolution by locally free sheaves.
If $C$ is arithmetically Cohen--Macaulay, then 
$F_3=0$ and $\sG$ splits into a direct sum of 
line bundles. In this case, the ideal $I_C$ is 
generated by the maximal minors of
$F_1 \leftarrow F_2$ by  the Hilbert--Burch 
Theorem~\cite{Hilb90,Burch,Eis}. In general, we 
have
\begin{align}
\label{eq:M_C props}
M_C \cong \sum_{n \in \ZZ} H^2(\PP^3,\sG(n))
\quad \text{and} \quad 
\sum_{n \in \ZZ}H^1(\PP^3,\sG(n))=0.
\end{align}
 
We explain now why curves linked by an even 
number of liaison steps have, up to a twist, the 
same Hartshorne--Rao module, thus illustrating 
its connection to liaison theory. We will then 
mention Rao's Theorem, which states that the 
converse also holds (Theorem~\ref{thm:rao}). 

Suppose that $f,g \in I_C$ are homogeneous forms  
of degree $d$ and $e$ without common factors. 
Let $X := V(f,g)$ denote the corresponding 
complete intersection, and let $C'$ be the 
residual scheme defined by the homogeneous ideal
$I_{C'}:=(f,g):I_C$ \cite{PS74}.
The locally free resolutions of $\sO_C$ and 
$\sO_{C'}$ are closely related, as follows. 
Applying $\mathcal Ext^2(-,\omega_{\PP^3})$ to 
the sequence
\[
0 \to \sI_{C/X} \to \sO_X \to \sO_C \to 0
\]
gives
\[
0 \leftarrow \mathcal Ext^2(\sI_{C/X}, 
\omega_{\PP^3}) \leftarrow \omega_X \leftarrow 
\omega_C \leftarrow 0.
\]
From $\omega_X \cong \sO_X(d+e-4)$, we conclude 
that 
$\mathcal Ext^2(\sI_{C/X},\sO_{\PP^3}(-d-e)) 
\cong \sO_{C'}$, and hence 
$\sI_{C'/X} \cong \omega_C(-d-e+4)$.
Now the mapping cone of 
\[
\xymatrix{
0 & \sO_C \ar[l] & \sO_{\PP^3} \ar[l] &
\bigoplus_j \sO_{\PP^3}(-j)^{\beta_{1j}}\ar[l] & \sG \ar[l]&0 \ar[l] \\
0 &\sO_X \ar[l] \ar[u]   &   \sO_{\PP^3} \ar[l] \ar[u]_\cong  & \sO_{\PP^3}(-d) \oplus \sO_{\PP^3}(-e) \ar[l] \ar[u] & \sO_{\PP^3}(-d-e) \ar[l] \ar[u]& 0  \ar[l] 
  }
\]
dualized with $\sHom(-,\sO_{\PP^3}(-d-e))$ gives
\[
\xymatrix{
0\to \sO_{\PP^3}(-d-e) \ar[r] \ar[d]_\cong&\bigoplus_j \sO_{\PP^3}(j-d-e)^{\beta_{1j}}\ar[r] \ar[d]& \sG^*(-d-e)\ar[d] \ar@{->>}[r]& \sI_{C'/X} \ar[d]   \\
0 \to   \sO_{\PP^3}(-d-e) \ar[r] & \sO_{\PP^3}(-e) \oplus \sO_{\PP^3}(-d) \ar[r] & \sO_{\PP^3} \ar@{->>}[r] & \sO_X, 
  }
\]
which yields the following locally free 
resolution of $\sO_{C'}$:
\[
0 \to 
\bigoplus_j \sO_{\PP^3}(j-d-e)^{\beta_{1j}} 
\to 
\sG^*(-d-e)\oplus \sO_{\PP^3}(-e) \oplus \sO_{\PP^3}(-d) 
\to  \sO_{\PP^3} \to  \sO_{C'} 
\to 0.
\]
In particular, after truncating this complex to  
resolve $I_{C'}$, one sees that 
\begin{align*}
M_{C'}&
:= \sum_{n\in \ZZ} H^1(\PP^3, \sI_{C'} (n)) 
\cong  \sum_{n \in \ZZ} H^1(\PP^3,\sG^*(n-d-e)) \\
& \phantom{:}\cong  \sum_{n\in \ZZ} H^2(\PP^3,\sG(d+e-4-n))^* \cong  \Hom_\bK(M_C,\bK)(4-d-e).
\end{align*}
Thus curves that are related via an even number 
of liaison steps have the same Hartshorne--Rao 
module up to a twist. Rao's famous result says 
that the converse is also true.

\begin{theorem}[Rao's Theorem~\cite{Rao78}]
\label{thm:rao} 
The even liaison classes of curves in $\PP^3$ 
are in bijection with finite length graded $S$-
modules up to twist.
\qed
\end{theorem}

Therefore the difficulty in constructing the 
desired space curve $C$ (of degree $d$ and genus 
$g$) lies completely in the construction of the 
appropriate Hartshorne--Rao module $M=M_C$.
Upon constructing $M$, we may then obtain the 
desired ideal sheaf $\sI_C$ as follows. 
Assume that we have a free $S$-resolution of 
$M_C$, 
\[ 
0 \leftarrow M_C \leftarrow F_0 \leftarrow F_1 \leftarrow F_2 \leftarrow F_3 \leftarrow F_4 \leftarrow 0,
\]
with $F_i= \bigoplus_j S(-j)^{\beta_{ij}}$. 
Let $\sF := \widetilde N$ be the sheafification  
of $N := \ker( F_1 \to F_0)$, the second syzygy 
module of $M$. In this case, $\sF$ will be a 
vector bundle without line bundle summands such 
that $H_*^1(\sF) \cong H^1_* (\sI_C)$ and 
$H_*^2(\sF)=0$. Here, we have used the notation 
$H_*^i(\sF) := \bigoplus_n H^i(\sF(n))$. If we 
constructed the correct Hartshorne--Rao module 
$M$, then taking $\sL_1$ and $\sL_2$ to be 
appropriate choices of direct sums of line 
bundles on $\PP^3$, a general homomorphism 
$\varphi \in \Hom(\sL_1,\sF\oplus \sL_2)$ will 
produce the desired curve $C$, as we will obtain 
$\sI_C$ as the cokernel of a map $\varphi$ of 
the bundles 
\[ 
\xymatrix{ 0 \ar[r] & \sL_1 \ar[r]^{\varphi\hspace{3.5mm}} & \sF \oplus \sL_2  \ar[r] &\sI_C \ar[r] & 0 \\
}.
\]

To compute the rank of $\sF$ and to choose the direct sums 
of line bundles $\sL_1$ and $\sL_2$, we now make 
plausible assumptions about the Hilbert function 
of $M_C$. We illustrate this approach in the 
example of the construction of a smooth linearly 
normal curve $C$ of degree $d=11$ and genus 
$g=10$. Since $ 2 d > 2g-2$, the line bundle 
$\sO_C(2)$ is already non-special. Hence by 
Riemann--Roch, we have that 
$h^0 (\sO_C(2))=22+1-10=13$.

\begin{remark}
\label{rem:maximal rank}
If we assume that $C$ is a curve of maximal 
rank, i.e., that all maps 
$H^0 (\sO_{\PP^3}(n)) \to H^0 (\sO_C(n))$ are 
either injective or surjective, then we can 
compute the Hilbert function of $M_C$ and $I_C$. 
Note that being of maximal rank is an open 
condition, so among the curves in the union 
$\sH_{d,g}$ of the component of the Hilbert 
scheme $\Hilb_{dt+1-g}(\PP^3)$  containing 
smooth curves, maximal rank curves form an open 
(and hopefully nonempty) subset.  There is a vast literature
on the existence of maximal rank curves; see, for example,~\cite{Flo}. 
\end{remark}

To gain insight into the Betti numbers of 
$M=M_C$, we use Hilbert's formula for the 
Hilbert series: 
\[
h_M(t) =  \sum_{n\in \ZZ} \dim M_n t^n
=\frac{ \sum_{i=0}^3 (-1)^i \sum_j \beta_{ij} t^j}{(1-t)^4}.
\]
Since $h_{M_C}(t)= 3t^2+4t^3$ by our maximal 
rank assumption 
(Remark~\ref{rem:maximal rank}), we have 
\begin{equation*}
(1-t)^4 h_{M}(t)= 
3t^2-8t^3+2t^4+12t^5-13t^6+4t^7, 
\end{equation*}
and thus the Betti table of $M$ must be 
\begin{center}
$\beta(M) \ = $ \ \ 
\begin{tabular}{c|c c c c c}
            & 0 &1 & 2 & 3 &4  \\ \hline
         2& 3 &8 &  2 &  .& .\\
         3&. &.& 12 &13& 4
\end{tabular}, 
\end{center}
\noindent
if we assume that $M$ has a so called  
\textbf{natural resolution}, which means that 
for each degree $j$ at most one $\beta_{ij}$ is 
nonzero. Note that having a natural resolution 
is an open condition in a family of modules with 
constant Hilbert function.

\begin{figure}
\begin{center}
\begin{tabular}{c||c|c|c|c}
$n$ &$h^1(\sI_C(n))$ & $h^0 (\sO_C(n))$ & $h^0 (\sO_{\PP^3}(n))$ & $h^0 (\sI_C(n))$   \cr \hline\hline
0 & 0 & 1 & 1 & 0 \cr
1 & 0 & 4 & 4 & 0 \cr
2 &  3  &  13  & 10 & 0 \cr
3 &   4 &   24 & 20 & 0 \cr
4 &   0   &   35 & 35 & 0 \cr
5 &  0 & 46 &  56 & 10 \cr
6 & 0 & 55 & 84 & 29 \cr
\end{tabular}
\end{center} 
\vspace{-5mm}
\caption{
With our maximal rank assumption of 
Remark~\ref{rem:maximal rank}, this table 
provides the relevant Hilbert functions in the 
case $d=11$ and $g=10$.}
\label{fig:HF M}
\end{figure}

Figure~\ref{fig:HF M} provides a detailed look 
at the Hilbert functions relevant to our 
computation. From these we see that 
$H^0_* (\sO_C)$ and $S_C=S/I_C$ will have the 
following potential Betti tables, if we assume 
that they also have natural resolutions: 

\begin{center}
$\beta(H^0_*(\cO_C)) \ =$ \ \ 
\begin{tabular}{c|c c c }
             &0 & 1& 2 \\ \hline
          0& 1& .& . \\
          1& . &. & . \\
          2& 3 & 8 & 2 \\
          3& .& . & 2\\
\end{tabular} 
$\quad$ and  $\quad \beta(S_C) \ =$ \ \ 
\begin{tabular}{c|c c c c }
             &0 & 1 &  2 & 3 \\ \hline
          0& 1  &. &  . & . \\
          1& . & .&  .& . \\
          2& . & . &  . & . \\
          3& .  & . & . & . \\
           4& . & 10 &13 & 4 \\
\end{tabular}. 
\end{center}

Comparing these Betti tables, we find the 
following plausible choices of $\sF$, $\sL_1$, 
and $\sL_2$:
\begin{itemize}
\vspace{-2mm}
\item 
We choose $\sF := \widetilde N$, where 
$N =\ker( \psi: S^8(-3) \to S^3(-2))$ is a 
sufficiently general $ 3\times 8$ matrix of 
linear forms; in particular, $\rank \sF =5$.

\vspace{-2mm}
\item 
Let $\sL_1:=\sO^2(-4)\oplus \sO^2(-5)$ and 
$\sL_2:=0$.
\vspace{-2mm}
\item 
The map $\varphi \in \Hom(\sL_1, \sF)$ is a 
sufficiently general homomorphism. Since the map 
$F_2 \to H^0_*(\sF)$ is surjective, the choice 
of $\varphi$ amounts to choosing an inclusion 
$\sO^2(-5) \to \sO^{12}(-5)$, i.e., a point in 
the Grassmannian $\GG(2,12)$.
\vspace{-2mm}
\item 
Finally, $\sI_C= \coker \varphi$.
\end{itemize}

It is not clear that general choices as above 
will necessarily yield a smooth curve. If the 
sheaf $\sHom(\sL_1,\sF \oplus \sL_2)$ happens to 
be generated by its global sections
$\Hom(\sL_1,\sF \oplus \sL_2)$, then a 
Bertini-type theorem as in~\cite{Kleiman} would 
apply. However, since we have to take all 
generators of $H^0_* (\sF)$ in degree $4$, this 
is not the case. On the other hand, there is no 
obvious reason that $\coker \varphi$ should not 
define a smooth curve, and upon construction, it 
is easy to check the smoothness of such an 
example using a computer algebra system, 
e.g., {\it Macaulay2} or {\sc Singular}. Doing 
this, we find that general choices do lead to a 
smooth curve.

\begin{exercise} 
Construct examples of curves of degree and genus 
as prescribed in Hartshone's book~\cite{Hart} in 
Figure~18 on page~354, including those which 
were open cases at the time of the book's 
publication.
\end{exercise}

\section{Random curves}
\label{sec:unirational}

In this section, we explain how the ideas of  
Section~\ref{sec:HRmods} lead to a 
computer-aided proof of the unirationality of 
the moduli space $\gM_g$ of curves of genus 
$g$, when $g$ is small. We will illustrate 
this approach by example, through the case of 
genus $g=12$ and degree $d=13$ in 
Theorem~\ref{thm:genus12,d13}.

\begin{definition} 
\label{def:unirational+}
A variety $X$ is called \textbf{unirational} if 
there exists a dominant rational map 
$\AA^n \dasharrow X$. A variety $X$ is called 
\textbf{uniruled} if there exists a dominant 
rational map $\AA^1\times Y \dasharrow X$ for 
some variety $Y$ that does not factor through 
$Y$. A smooth projective variety $X$ has 
\textbf{Kodaira dimension} $\kappa$  if the 
section ring 
$R_X:=\sum_{n\ge 0} H^0(X,\omega_X^{\tensor n})$ 
of pluri-canonical forms on $X$ has a Hilbert 
function with growth rate 
$h^0(\omega_X^{\tensor n}) \in O(n^\kappa)$. We 
say that $X$ has \textbf{general type} if 
$\kappa = \dim X$, the maximal possible value. 
\end{definition}

Since the pluri-genera 
$h^0(\omega_X^{\tensor n})$ are birational 
invariants, being of general type does not 
depend on a choice of a smooth compactification. 
Thus we may also speak of general type for 
quasi-projective varieties.

Unirationality and general type are on opposite 
ends of birational geometry. If a variety is of 
general type, then there exists no rational 
curve through a general point of 
$X$~\cite[Corollary~IV.1.11]{Kol}.
On the other hand, uniruled varieties have the 
pluri-canonical ring $R_X=(R_X)_0=\CC$ and thus 
(by convention) have Kodaira dimension 
$\kappa=-\infty$. In fact, even if $X$ is 
unirational, then we can connect any two general 
points of $X$ by a rational curve. 

We now recall results concerning the 
unirationality of the moduli space $\gM_g$. 
There are positive results for small genus, 
followed by negative results for large genus. 

\begin{theorem}[Severi, Sernesi, Chang--Ran, Verra]
\label{thm:g14} 
The moduli space $\gM_g$ of curves of genus $g$ 
is unirational for $g\le 14$.
(For $g\le 10$, see~\cite{Sev}. 
For $g=12,11,13$, see~\cite{Ser,CR}. For $g=14$, 
see~\cite{Ve}.)
\qed
\end{theorem}

\begin{theorem}[Harris--Mumford, Eisenbud--Harris, Farkas~\cite{HM,EH,Fa,Fa1,Fa2}]
\label{thm:22,24+}
The moduli space $\gM_g$ of curves of genus $g$ 
is of general type for $g\ge 24$ or $g=22$. The 
moduli space $M_{23}$ has  Kodaira 
dimension $\ge 2$.
\qed
\end{theorem}

We call this beautiful theorem a negative result 
because it says that it will be very difficult 
to write down explicitly a general curve of 
large genus. Given a family of curves of genus 
$g\ge 24$ that pass through a general point of 
$\gM_g$, say via an explicit system of equations  
with varying coefficients, none of the  
essential coefficients is a free parameter. All 
of the coefficients will satisfy some 
complicated algebraic relations.
On the other hand, in unirational cases, there 
exists  a dominant family of curves whose 
parameters  vary freely.

In principle, we can compute a dominating family 
explicitly along with a unirationality proof. In 
practice, this is often out of reach using 
current computer algebra systems; however, the 
following approach is feasible today in many 
cases. 
By replacing each free parameter in the 
construction of the family by a randomly chosen 
value in the ground field, the computation of an 
explicit example is possible. 
In particular, over a finite field $\FF$, where 
it is natural to use the constant probability 
distribution on $\FF$, a unirationality proof 
brings with it the possibility of choosing 
random points in $\gM_g(\FF)$, i.e., to compute 
a random curve. These curves can then be used 
for further investigations of the moduli space, 
as well as to considerably simplify the existing 
unirationality proofs. The advantage of using 
such random curves in the unirationality proof 
is that, with high probability, they will be 
smooth curves, while in a theoretical treatment, 
smoothness is always a delicate issue. 

To begin this construction, we first need some 
information on the  projective models of a 
general curve. This is the content of Brill--
Noether theory. Let
\[ 
W^r_d(C) := 
\{ L \in \Pic^d(C) \mid h^0(C,L) \ge r+1 \} 
\subset \Pic^d(C)
\]
denote the space of line bundles of degree $d$ on $C$ 
that give rise to a morphism $C \to \PP^r$.

\begin{theorem}[Brill--Noether, Griffith--Harris, Fulton--Lazarsfeld, Gieseker] \label{thm:BNGHFGL}
Let $C$ be a smooth projective curve of genus 
$g$. 
\begin{enumerate}
\vspace{-1mm}
\item 
~\emph{\cite{BN}}
At every point, 
$\dim W^r_d(C) \ge \rho := g-(r+1)(g-d+r)$. 
\vspace{-1mm}
\item 
~\emph{\cite{GH,Ful}}
If $\rho \ge 0$, then $W^r_d(C) \not=0$, and if 
$\rho >0$, then $W^r_d(C)$ is connected. 
Further, the tangent space of $W^r_d(C)$ at a 
point $L \in W^r_d(C) \setminus W^{r+1}_d(C)$ is 
\[
T_L W^r_d(C) = 
{\rm Im}\  \mu_L^\perp \subset H^1(\sO_C) 
= T_L \Pic^d(C),
\] 
where 
$\mu_L\colon H^0(L)\tensor 
H^0(\omega_C\tensor L^{-1}) 
\to H^0(\omega_C) = H^1(\sO_C)^*$ 
denotes the Petri map.
\vspace{-1mm}
\item 
~\emph{\cite{Gie}}
If $C \in \gM_g$ is a general curve, then 
$W^r_d(C)$ is smooth of dimension $\rho$ away 
from $W^{r+1}_d(C)$. More precisely, the Petri 
map $\mu_L$ is injective for all 
$L \in W^r_d(C)\setminus W^{r+1}_d(C)$.
\qed
\end{enumerate}
\end{theorem}

We now illustrate the computer-aided 
unirationality proof of $\gM_g$ by example, 
through the case $g=12,d=13$ \cite{ST}. This case is not 
amongst those covered by Sernesi~\cite{Ser} or 
Chang--Ran~\cite{CR}. They treated the cases 
$g=11,d=12$, $g=12,d=12$, and $g=13,d=13$. We 
are choosing the case $d=12,g=13$ because it 
illustrates well the difficulty of this 
construction.  For $g=14$, see~\cite{Ve} and, 
for a computer aided unirationality 
proof,~\cite{Sch4}. For a related {\it Macaulay2} package, 
see~\cite{BGS}. 

\begin{theorem}
\label{thm:genus12,d13}
Let $g=12$ and $d=13$. Then 
$\Hilb_{dt+1-g}(\PP^3)$ has a component 
$\sH_{d,g}$ that is unirational and dominates 
the moduli space $\gM_g$ of curves of genus $g$.
\end{theorem}

\begin{proof}
This proof proceeds as follows. We first compute 
the Hilbert function and expected syzygies of 
the Hartshorne--Rao module $M=H^1_*(\sI_C)$, the 
coordinate ring $S_C$, and the section ring 
$R:=H^0_*(\sO_C)$. We then use this information 
to choose generic matrices which realize the 
free resolution of $M$. Finally, we show that 
this construction leads to a family of curves 
that dominate $\gM_{12}$ and generically 
contains smooth curves.

We first choose $r$ so that a general curve has 
a model of degree $d=13$ in $\PP^r$. In our 
case, we choose $r=3$ so that $g-d+r=2$.
To compute the Hilbert function and expected 
syzygies of the Hartshorne--Rao module 
$M=H^1_*(\sI_C)$, the coordinate ring 
$S_C=S/I_C$, and the section ring 
$R=H^0_*(\sO_C)$, we assume the open condition 
that $C$ has maximal rank, i.e.,
\[
H^0(\PP^3,\sO(n)) \to H^0(C, L^n)
\]
is of maximal rank for all $n$, as in 
Remark~\ref{rem:maximal rank}. 
In this case, $h_M(t)=5t^2+8t^3+6t^4$, which has 
Hilbert numerator 
\[
h_M(t)(1-t)^4
= 5t^2-12t^3+4t^4+4t^5+9t^6-16t^{10}+6t^{11}. 
\]
If $M$ has a natural resolution, so that for 
each $j$ at most one $\beta_{ij}(M)$ is nonzero, 
then $M$ has the Betti table 
\begin{center}
$\beta(M) \ = $ \ \ 
\begin{tabular}{c|ccccc}
           & 0 &  1 & 2 & 3 &4 \\ \hline
         2 & 5 &12 & 4 &  . & . \\
         3& . & .  &4  &. &.\\
         4& . & . & 9 &16& 6\\
\end{tabular}.
\end{center}
\noindent
If we assume the open condition that $S_C$ and 
$R$ have natural syzygies as well, then their 
Betti tables are
\begin{center}
$\beta(S_C) \ = $ \ \ 
\begin{tabular}{c|c c c }
             &0 & 1& 2 \\ \hline
          0& 1& .& . \\
          1& . &. & . \\
          2& 3 & 12 & 4 \\
          3& .& . & 2\\
\end{tabular} $\quad$ and  $\quad$
$\beta(R) \ = $ \ \ 
\begin{tabular}{c|c c c c }
             &0 & 1 &  2 & 3 \\ \hline
          0& 1  &. &  . & . \\
          1& . & .&  .& . \\
          2& . & . &  . & . \\
          3& .  & . & . & . \\
           4& . & 2 &. & . \\
           5& . & 9 & 16 & 6 \\
\end{tabular}. 
\end{center}
We conclude that once we have constructed the 
Hartshorne--Rao module $M=M_C$, say via its 
representation
\[ 
0 \leftarrow M  \leftarrow S^5(-2) 
\leftarrow S^{12}(-3),
\]
we may choose $\sF$ to be the kernel of
\[ 
0 \leftarrow \sO^5(-2) \leftarrow \sO^{12}(-3) 
\leftarrow \sF \leftarrow 0
\]
and set $\sL_1 := \sO(-4)^4 \oplus \sO^2(-5)$ 
and $\sL_2:=0$.
Then $C$ is determined by $M$ and the choice of 
a point in $\GG(2,4)$. In particular, as 
mentioned earlier, constructing $C$ is 
equivalent to constructing the finite length 
module $M$ with the desired syzygies.

If we choose the presentation matrix $\phi$ of 
$M$ to be given by a general (or random) 
$5\times 12$ matrix of linear forms, then its 
cokernel will be a module with Hilbert series 
$5t^2+8t^3+2t^4$. In other words, to get the 
right Hilbert function for $M$, we must force 
$4$ linear syzygies. 
To do this, choose a general (or random) 
$12\times 4$ matrix $\psi$ of linear forms. Then 
\[
\ker (\psi^t\colon S^{12}(1) \to S^4(2)) 
\] 
has at least $12\cdot 4- 4 \cdot 10 =8$ 
generators in degree $0$.
In fact, there are precisely $8$ and a general 
point in $\GG(5,8)$ gives rise to a 
$12 \times 5$ matrix $\varphi^t$ of linear 
forms. This means that 
$M := \coker(\varphi\colon S^12(-3)\to S^5(-2))$
to have Hilbert series $5t^2+8t^3+6t^4$, due to 
the forced $4$ linear syzygies.

Having constructed $M$, it remains to prove that 
this construction leads to a family of curves 
that dominates $\gM_{12}$. 
To this end, we compute a random example $C$, 
say over a finite prime field $\FF_p$, and 
confirm its smoothness. Since we may regard our 
computation over $\FF_p$ as the reduction modulo 
$p$ of a construction defined over an open part 
of $\Spec\ZZ$, semi-continuity allows us to 
establish the existence of a smooth example 
defined over $\QQ$ with the same syzygies. 

We now consider the universal family 
$\gW^r_d \subset \gP ic^d$ over $\gM_g$ and a 
neighborhood of our example 
$(C,L) \in \gP ic^d$. Note that the codimension 
of $\gW^r_d$ is at most 
$(r+1)(g-d+r)=4\cdot 2 =8$. 
On the other hand, we claim that the Petri map 
$\mu_L$ for ($C,L)$ is injective. (Recall the 
definition of $\mu_L$ from 
Theorem~\ref{thm:BNGHFGL}.) To see this, note 
that the Betti numbers of $H^0_*(\omega_C)$ 
correspond to the dual of the resolution of 
$H^0_*(\sO_C)$, so 
\begin{center}
$\beta(H^0_*(\omega_C)) \ = $ \ \ 
\begin{tabular}{c|c c c }
             &0 & 1& 2 \\ \hline
          -1& 2& .& . \\
          0& 4& 12 & 3 \\
          1& . & . & . \\
          2& .& . & 1\\
\end{tabular}.
\end{center}
Thus there are no linear relations among the two 
generators in $H^0(\omega_C \tensor L^{-1})$, 
which means that the $\mu_L\colon 
H^0(L) \tensor H^0(\omega_C \tensor L^{-1})$ 
is injective. From this we see that 
$\dim W^r_d(C)$ has dimension $4$ at $(C,L)$, 
and the constructed family dominates for 
dimension reasons. 
\end{proof}

The unirationality of $\gM_{15}$ and $\gM_{16}$ 
are open; however, these moduli spaces are 
uniruled.

\begin{theorem}[Chang--Ran~\cite{CR1,CR2,BV,Fa1}]
\label{thm:15,16uniruled}
The moduli space $\gM_{15}$ is rationally connected, and $\gM_{16}$ is uniruled.
\qed
\end{theorem}

To explain why the unirationality in these cases 
is more difficult to approach using the method 
of Theorem~\ref{thm:genus12,d13}, we conclude 
with a brief discussion on the space models of 
curves of genus $g=16$. By Brill--Noether 
theory, a general curve $C$ of genus $16$ has 
finitely many models of  degree $d=15$ in 
$\PP^3$. Again assuming the maximal rank 
condition of Remark~\ref{rem:maximal rank}, the 
Hartshorne--Rao module $M=H^1_*(\sI_C)$ has 
Hilbert series
\[
H_M(t)= 5t^2+10t^3+10t^4+4t^5
\]
and expected syzygies
\begin{center}
$\beta(M) \ = $ \ \ 
\begin{tabular}{c|c c c cc}
             &0 & 1& 2 & 3 &4\\ \hline
          2& 5& 10& . \\
          3& .& . & 4 &.&.\\
          4& . & . & 9 & 6 &. \\
          5& .& . &  . & 6 & 4\\
\end{tabular}.
\end{center}
The section ring $H^0_*(\sO_C)$ and the 
coordinate ring $S_C$ have expected syzygies
\begin{center}
$\beta(H^0 _*(\cO_C)) \ = $ \ \ 
\begin{tabular}{c|c c c }
             &0 & 1& 2 \\ \hline
          0& 1& .& . \\
          1& . &. & . \\
          2& 5 & 10& . \\
          3& .& . & 4\\
\end{tabular} $\quad$ and  $\quad$
$\beta(S_C) \ = $ \ \ 
\begin{tabular}{c|c c c c }
             &0 & 1 &  2 & 3 \\ \hline
          0& 1  &. &  . & . \\
          1& . & .&  .& . \\
          2& . & . &  . & . \\
          3& .  & . & . & . \\
           4& . & . &. & . \\
           5& . & 9 & 6 & . \\
           6& . & . & 6 & 4
\end{tabular}. 
\end{center}

\begin{proposition}
\label{prop:genus16} 
A general curve $C$ of genus $g=16$ and degree 
$d=15$ in $\PP^3$ has syzygies as above. In 
particular, the Hartshorne--Rao module $M_C$ 
uniquely determines $C$. 
Furthermore,  the rational map from the 
component $\sH_{d,g}$ of the Hilbert scheme 
$\Hilb_{15t+1-16}(\PP^3)$ that dominates 
$\gM_{16}$ defined by
\begin{align*}
\sH_{d,g} &
	\dasharrow \{ \hbox{ 20 determinantal 
	points } \}\\ 
C &
	\,\,\mapsto \,\, 
	\Gamma:=\supp \coker( \varphi^t: \sO^6(-1) 
	\to \sO^4)  
\end{align*}
is dominant. Here 
$\varphi\colon S^4(-9) \to S^6(-8)$ denotes the 
linear part of the last syzygy matrix of $M$.
\end{proposition}
\begin{proof} 
For the first statement, it suffices to find an 
example with the expected syzygies, since Betti 
numbers behave semi-continuously in a family of 
modules with constant Hilbert function. We may 
even take a reducible example, provided that it 
is smoothable. Consider the union 
$C:=E_1 \cup E_2 \cup E_3$ of three smooth 
curves of genus $2$ and degree $5$, such that 
$E_i \cap E_j$ for $i \not=j$ consists of $4$ 
nodes of $C$.  Then $C$ has degree 
$d=3\cdot 5=15$ and genus 
$g=3\cdot 2 + 4\cdot 3-2=16$. Clearly, 
$C$ is smoothable as an abstract curve. For 
general choices, it is smoothable as an embedded 
curve because the $g^3_{15}$ on the reducible 
curve is an isolated smooth point in $W^3_{15}$ 
(as we will see), so the smooth curves nearby 
have an isolated $g^3_{15}$ as well.

It is easy to find such a union over a finite 
field $\FF$. Start with the $12$ intersection 
points $\{p_1,\ldots,p_4\} \cup 
\{p_5,\ldots,p_8\} \cup \{p_9,\ldots,p_{12} \}$ 
randomly chosen in $\PP^3(\FF)$. To construct 
$E_1$, pick at random a quadric $Q_1$ in the 
pencil of quadrics through $\{p_1,\ldots,p_8\}$.  
Next, we must check if the tangent hyperplane of 
$Q_1$ in a point, say $p_1$, intersects $Q_1$ in 
a pair of lines individually defined over $\FF$; 
this will happen about $50 \%$ of the time. Once 
this is true, choose one of the lines, call it 
$L_1$. Then 
$| \sO_{Q_1}(3) \tensor \sO_{Q_1}(-L_1)|$
is a linear system of class $(3,2)$ on 
$Q_1 \cong \PP^1 \times \PP^1$. We may take 
$E_1$ as a general curve in this linear system 
that passes through  $\{p_1,\ldots,p_8\}$.
Similarly, we choose $E_2$ using 
$\{p_1,\ldots,p_4,p_9,\ldots,p_{12} \}$ and 
$E_3$ starting with $\{p_5, \ldots,p_{12}\}$. 
The union of the $E_i$ yields the desired curve 
$C$, and it a straightforward computation to 
check that $C$ has the expected Hartshorne--Rao 
module and syzygies. 

The second statement can be proved by showing 
that the appropriate map between tangent spaces 
is surjective for this example. This involves 
computing appropriate $\Ext$-groups. Define 
\begin{align*}
\overline M :&= \coker( S^6(-2)\oplus S^6(-1) 
\to S^4)\\
& =\Ext^4_S(M,S(-9))\\
&= \Hom_K(M,K)(-5)\\
\text{and}\ \ \, 
N:&=\coker( \varphi^t\colon S^6(-1)\to S^4).
\end{align*}
Then there is a short exact sequence
\[ 
0 \to P \to N \to \overline M \to 0
\]  
of modules with Hilbert series 
\begin{eqnarray*}
h_{\overline M}(t) &=&4 +10t+10t^2 +5t^3\\
h_N(t) &=& 4+10t+16t^2+20t^3+20t^4+20t^5+\cdots  \\
h_P(t)& = & \qquad \qquad \; \; 
6t^2+15t^3+20t^4+20t^5+\cdots.
\end{eqnarray*}
The group $\Ext^1_S(\overline M, \overline M)$ 
governs the deformation theory of $\overline M$ 
(and $M$). More details can be found 
in~\cite{Hart-def}, for example, 
Theorem~2.7 applied in the affine case.
More precisely, the degree $0$ part of this 
$\Ext$-group is the tangent space of homogeneous  
deformations of $M$, which in turn is isomorphic 
to the tangent space of the Hilbert scheme in 
$C$. Similarly, in the given example, 
$\Ext^1(N,N)_0$ can be identified with the 
tangent space to the space of twenty 
determinantal points. Note that we have the 
diagram
\[  
\xymatrix{
        \Ext^1_S(\overline M, \overline M) \ar[r]             &   \Ext^1_S(N,\overline M) \ar[r] & \Ext^1_S(P,\overline M). \\
   & \Ext^1_S(N,N) \ar[u] & \\
   & \Ext^1_S(N,P) \ar[u] &\\  }
\]
In our example, computation shows that 
\begin{align*}
\label{eq:Exts}
\nonumber
\dim\Ext^1_S(\overline M,\overline M)_0&=60, \\
\dim\Ext^1_S(N,\overline M)_0 
= \dim  \Ext^1_S(N,N)_0 &= 45, 
\quad \text{and}\\
\dim \Ext^1_S(P,\overline M)_0
=\dim \Ext^1_S(N,P)_0&=0. 
\end{align*} 
Thus the induced map 
$\Ext^1_S(\overline M, \overline M)_0 
\to \dim  \Ext^1_S(N,N)_0$ is surjective with 
$15$-dimen\-sional kernel, as expected.
\end{proof}

\begin{exercise} Fill in the computational 
details in of the proof of 
Proposition~\ref{prop:genus16} and 
Theorem~\ref{thm:genus12,d13} using your 
favorite computer algebra system, say
{\it  Macaulay2} or {\sc Singular}. 
\end{exercise}

\begin{remark} 
In the proof of Proposition~\ref{prop:genus16}, 
the module $P$ has syzygies
\begin{center}
$\beta(P) \ = $ \ \ 
\begin{tabular}{c|c c c cc}
             &0 & 1& 2 & 3 \\ \hline
          2& 6& 9& . &.\\
          3& .& 4 & 6 &.\\
          4& . & . & 6 & 5 \\
\end{tabular}.
\end{center}
The cokernel of $\psi^t\colon S^6(-1) \to S^5$ 
has support on a determinantal curve $E$ of 
degree $15$ and genus $26$, which is smooth for 
general $C$. The points $\Gamma$ form a divisor 
on $E$ with $h^0(E,\sO_E(\Gamma))=1$. The curves 
$E$ and $C$ do not intersect; in fact, we have 
no idea how the curve $E$ is related to $C$, 
other than the fact that it can be constructed 
from the syzygies of $M$. It is possible that 
$\gM_{16}$ is not unirational, and, even if 
$\gM_{16}$ is unirational, it could be that the 
component of the Hilbert scheme containing $C$ 
is itself not unirational. 

It is not clear to us whether it is a good idea 
to start with the determinantal points $\Gamma$ 
in Proposition~\ref{prop:genus16}. Perhaps 
entirely different purely algebraic methods 
might lead to a unirational construction of the 
modules $M$, and we invite the reader to 
discover such an approach. 
\end{remark}


\bigskip

\vbox{\noindent Author Addresses:\par
\medskip

\noindent{Christine Berkesch}\par
\noindent{Department of Mathematics, Duke University, Box 90320\\ Durham, NC 27708}\par
\noindent{cberkesc@math.duke.edu}\par

\medskip
\noindent{Frank-Olaf Schreyer}\par
\noindent{Mathematik und Informatik, Universit\"at des Saarlandes, Campus E2 4, 
D-66123 Saarbr\"ucken, Germany}\par
\noindent{schreyer@math.uni-sb.de}\par
}
\end{document}

%% file: preamble.tex
\def\antiddot{\mathinner{\mkern1mu\raise1pt\vbox{\kern7pt\hbox{.}}\mkern2mu
        \raise4pt\hbox{.}\mkern2mu\raise7pt\hbox{.}\mkern1mu}}

\newcommand{\bK}{{\Bbbk}}

\newcommand{\CC}{{\mathbb C}}

\newcommand{\FF}{{\mathbb F}}
\newcommand{\GG}{{\mathbb G}}

\newcommand{\NN}{{\mathbb N}}
\newcommand{\PP}{{\mathbb P}}
\newcommand{\QQ}{{\mathbb Q}}

\newcommand{\ZZ}{{\mathbb Z}}

\newcommand{\Ext}{{\rm{Ext}}}

\newcommand{\ini}{{\bf{L}}}

\newcommand{\fr}{{f_1,\ldots,f_r}}
\newcommand{\PR}{{\bK [x_1, \dots , x_n]}}

\newcommand{\coker}{{\rm{coker}\,}}

\newcommand{\s}{\mathcal}

\newcommand{\sF}{{\s F}}
\newcommand{\sG}{{\s G}}
\newcommand{\sH}{{\s H}}
\newcommand{\sI}{{\s I}}

\newcommand{\sL}{{\s L}}
\newcommand{\sM}{{\s M}}

\newcommand{\sO}{{\s O}}

\newcommand{\cO}{{\s O}}

\newcommand{\lra}{\longrightarrow}

\newcommand{\tensor}{\otimes}

\newcommand{\punkt}{\hspace{-.3ex}\raise.15ex\hbox to1ex{\Huge.}}

\DeclareMathOperator{\Pic}{Pic}

\DeclareMathOperator{\Spec}{Spec}

\DeclareMathOperator{\Hom}{Hom}
\DeclareMathOperator{\sHom}{\sH om}

\DeclareMathOperator{\supp}{supp}

\DeclareMathOperator{\im}{im}

\DeclareMathOperator{\pd}{pd}

\DeclareMathOperator{\codim}{codim}
\DeclareMathOperator{\rank}{rank}

\newcommand{\gP}{\mathfrak P}

\newtheorem{theorem}{Theorem}[section]

\newtheorem{proposition}[theorem]{Proposition}
\newtheorem{corollary}[theorem]{Corollary}
\newtheorem{conjecture}{Conjecture}[section]
\theoremstyle{definition}
\newtheorem{definition}[theorem]{Definition}

\newtheorem{remark}[theorem]{Remark}

\newtheorem{example}[theorem]{Example}
\newtheorem{exercise}[theorem]{Exercise}